
\baselineskip=14pt
\parskip=10pt

\font\eightrm=cmr8 

\magnification=\magstephalf
\def\W{{\cal W}}

\def\1{{\overline{1}}}
\def\2{{\overline{2}}}
\parindent=0pt
\overfullrule=0in

\def\frac#1#2{{#1 \over #2}}

\centerline
{
 \bf A Treatise on Sucker's Bets
}
\rm
\bigskip
\centerline
{\it By Shalosh B. EKHAD and Doron ZEILBERGER}

{\bf Preface}

Consider a gambling game where each player has a  die marked with dots (not necessarily the usual way).
For each die, each face is equally likely to show up, in other words the dice are {\it fair}.
Each player rolls his  die simultaneously, and whoever has  more dots on his landed face wins the round.

In 1970, Statistics giant, Bradley Efron, amazed the world ([G], ch. 22; [W]) by coming up with a set of four dice, 
let's call them $A$,$B$,$C$,$D$, whose faces are marked as follows
$$
A=[1,1,5,5,5,5] \quad , \quad
B=[4,4,4,4,4,4] \quad , \quad
C=[3,3,3,3,7,7] \quad , \quad
D=[2,2,2,6,6,6] \quad .
$$

It turns out that 

$\bullet$ die $A$ beats die $B$ (in the sense that $A$'s chance of winning exceeds $B$'s chance of winning) ;

$\bullet$ die $B$ beats die $C$ ;

$\bullet$  die $C$ beats die $D$ ; 

but,  {\bf surprise surprise}, 

$\bullet$ die D beats die A! 

This was an amazing demonstration that ``being more likely to win" is {\bf not} a transitive relation. 
But that was only one example, and of course, instead of dice, we can use decks of cards, 
that Martin Gardner ([G], ch. 23) called {\bf sucker's bets}.

Gardner gave an example of a set of three decks of three cards each, all marked with different numbers, that
is derived from the rows of a three by three magic square
$$
A=[1,6,8] \quad , \quad
B=[3,5,7] \quad , \quad
C=[2,4,9] \quad .
$$

Can you find {\bf all} such examples, with a specified number of decks, and deck-sizes? If you have a computer algebra system (in our case Maple), you sure can!

{\it Not only that}, we can figure out how likely such sucker's bets are, and derive, {\bf fully automatically}, statistical information! 

\vfill\eject

{\bf Why are they called Sucker's Bets}

The {\it sucker} tacitly assumes that `being more likely to win'  is a transitive relationship, hence he or she would not object
to the privilege of being {\bf first} to pick which deck to play with. Then the  {\it hustler} can always pick a better deck.

$\bullet$ If the sucker picks $A$, then the hustler will pick $C$ ;

$\bullet$ If the sucker picks $B$, then the hustler will pick $A$ ;

$\bullet$ If the sucker picks $C$, then the hustler will pick $B$ .

Let's first convince ourselves that $A$ is better than $B$, $B$ is better than $C$, and $C$ is better than $A$.

{\bf $A$ versus $B$}

Out of the {\bf nine} (equally likely!) possibilities $[a,b]$
$$
\{ [1,3],[1,5],[1,7] ,[6,3],[6,5],[6,7],[8,3],[8,4],[8,7] \} \quad,
$$
$A$ wins in $5$ of them, namely $\{[6,3],[6,5],[8,3],[8,4],[8,7] \}$,
while $B$ wins in $4$ of them, namely $\{ [1,3],[1,5],[1,7] , [6,7] \}$.

{\bf $B$ versus $C$}

Out of the {\bf nine} (equally likely!) possibilities $[b,c]$
$$
\{ [3,2],[3,4],[3,9] ,[5,2],[5,4],[5,9],[7,2],[7,4],[7,9] \} \quad,
$$
$B$ wins in $5$ of them, namely $\{[3,2],[5,2],[5,4],[7,2],[7,4] \}$,
while $C$ wins in $4$ of them, namely $\{ [3,4],[3,9],[5,9] , [7,9] \}$.

{\bf $C$ versus $A$}

Out of the {\bf nine} (equally likely!) possibilities $[c,a]$
$$
\{ [2,1],[2,6],[2,8] ,[4,1],[4,6],[4,8],[9,1],[9,6],[9,8] \} \quad,
$$
$C$ wins in $5$ of them, namely $\{[2,1],[4,1],[9,1],[9,6],[9,8] \}$,
while $A$ wins in $4$ of them, namely $\{ [2,6],[2,8],[4,6] , [4,8] \}$.

{\bf Other Examples of Surprising Non-Transitive Relations}

Of course, non-transitivity is nothing new! An ancient example is {\it love}. 
Joe loves his wife, Joe's wife loves her mother, but Joe {\bf hates} his mother-in-law.

Another example is the relation {\it being afraid of}. Indeed

$\bullet$ My dog is afraid of me ;

$\bullet$ My cat is afraid of my dog ;

$\bullet$  The mouse  is afraid of my cat ;

$\bullet$  My wife  is afraid of the mouse ;

$\bullet$  I am afraid of my wife.

Another, more ``serious'' example is in voting, where famously  the 18th-century French polymath,

{\it Marie Jean Antoine Nicolas de Caritat, Marquis de Condorcet}, 

came up with his famous {\bf paradox}, also mentioned by Gardner([G], ch. 23), 
that inspired Kenneth Arrow's {\it impossibility theorem}  (that earned him a Nobel!).

Another example is the intriguing {\bf Penny ante} ([P], see also [G], ch. 23, and  discussion in [NZ]).

{\bf  Mapping Sucker's Bets to Words}

First let's consider the case where there are $k$ decks, where deck $1$ has $a_1$ cards,
deck $2$ has $a_2$ cards, $\dots$, deck $k$ has $a_k$ cards, and let $N=a_1+ \dots +a_k$ be the
total number of cards participating. Let's first treat the case where all the denominations of
these $N$ cards are different, so without loss of generality, we can make them $\{1, 2, \dots , N\}$.
There are $(a_1+ \dots +a_k)!/(a_1! \cdots a_k!)$ ways of assigning the cards to the various decks,
and there is an obvious bijection between such decks and words $w$ in the alphabet $\{1, 2, \dots , k\}$
with $a_1$ occurrences of $1$, $a_2$ occurrences of $2$, $\dots$  $a_k$ occurrences of $k$,
where $w_i=j$ means that we put the card with the denomination $i$ into the $j$-th deck.

For pedagogical clarity, until further notice, let's take $k=3$.

For example, the above-mentioned set of three decks, each with three cards
$$
A=[1,6,8] \quad , \quad
B=[3,5,7] \quad , \quad
C=[2,4,9] \quad ,
$$
corresponds to the `word'
$$
ACBCBABAC \quad,
$$
and replacing $A,B,C$ by $1,2,3$, respectively, we get the `word'
$$
132321213 \quad .
$$

{\bf Which Words Correspond to Sucker's Bets?}

Let's consider the "magic" deck above, whose `word' turned out to be  $132321213$. The reason that it corresponds to a sucker's bet
is that 

$\bullet$ The number of times  letter `$1$' is to the left of letter `$2$' (not necessarily immediately before) 
is less than the number of times that a `$2$' is ahead of a `$1$'.

{\eightrm
(Indeed there are three `$2$'s after the `$1$' at the first place, and one `$2$' after the `$1$' at the sixth place, totaling
four occurrences of `$1$' before `$2$', while there are two `$1$'s that occur after the `$2$' at the third place,
two `$1$'s that occur after the `$2$' at the fifth place, and one `$1$' that occurs after the `$2$' at the seventh place.)}

$\bullet$ The number of times  letter `$2$' is to the left of letter `$3$' (not necessarily immediately before) 
is less than the number of times that a `$3$' is ahead of a `$2$'. (Check!)

$\bullet$ The number of times letter `$3$' is to the left of letter `$1$' (not necessarily immediately before) 
is less than the number of times that a `$1$' is ahead of a `$3$'. (Check!)

This leads us to introduce three {\bf word statistics}, for any word in the alphabet $\{1,2,3\}$
(below $|S|$ means, as usual, the number of elements in the set $S$).
$$
s_1(w) := |\{ (i,j) \, | \, i<j, w_i=2 \quad AND \quad w_j=1\}| \, - \,|\{ (i,j) \, | \, i<j, w_i=1 \quad AND \quad w_j=2\}|  \quad ,
$$
$$
s_2(w) := |\{ (i,j) \, | \, i<j, w_i=3 \quad AND \quad w_j=2\}| \, - \,|\{ (i,j) \, | \, i<j, w_i=2 \quad AND \quad w_j=3\}| \quad ,
$$
$$
s_3(w) := |\{ (i,j) \, | \, i<j, w_i=1 \quad AND \quad w_j=3\}| \, - \,|\{ (i,j) \, | \, i<j, w_i=3 \quad AND \quad w_j=1\}|  \quad .
$$
So for the above-mentioned $w=132321213$, $s_1(w)=s_2(w)=s_3(w)=1$.

This leads us to a {\bf characterization} of words that correspond to sucker's bets.

{\bf Proposition}: A word $w$ in the alphabet $\{1,2,3\}$ corresponds to a sucker's bets
(where the first deck has denominations indicating the locations of the letter $1$,
the second deck has denominations indicating the locations of the letter $2$,
and the third deck has denominations indicating the locations of the letter $3$)
{\bf if and only if}
$$
s_1(w) >0 \quad, \quad s_2(w) >0 \quad, \quad s_3(w)>0  \quad .
\eqno(SBC)
$$

{\bf Counting Sucker's Bets}

Given three positive integers $a_1,a_2,a_3$, out of the $(a_1+a_2+a_3)!/(a_1!a_2!a_3!)$ words $w$ in the
alphabet $\{1,2,3\}$ with $a_1$ $1$s,  $a_2$ $2$s,  $a_3$ $3$, how many have the property $(SBC)$?

In other words:

{\it
How many three-deck sets of cards form a sucker's bet where
the first deck has $a_1$ cards, the
second deck has $a_2$ cards, and the third deck has $a_3$ cards, 
and all denominations are different and are drawn from $\{1, 2, \dots, a_1+a_2+a_3\}$?
}

The naive way would be to actually examine all the  $(a_1+a_2+a_3)!/(a_1!a_2!a_3!)$ possible words, 
compute $s_1(w),s_2(w),s_3(w)$, and count those for which condition $(SBC)$ holds.
But there is a better way!

Let $q_1,q_2,q_3$ be three (commuting) {\it indeterminates}, and define the {\bf weight} of a word $w$ by
$$
weight(w):=q_1^{s_1(w)} q_2^{s_2(w)} q_3^{s_3(w)} \quad .
\eqno(Weight)
$$

Let $\W(a_1,a_2,a_3)$ be the set of words in $\{1,2,3\}$ with $a_1$ $1$s,$a_2$ $2$s, and $a_3$ $3$s.

Now define the {\it weight enumerator} (aka {\it generating function}) 
$$
F(a_1,a_2,a_3)(q_1,q_2,q_3) \, := \, \sum_{w \in \W(a_1,a_2,a_3)} weight(w) \, = \,
\sum_{w \in \W(a_1,a_2,a_3)} q_1^{s_1(w)} q_2^{s_2(w)} q_3^{s_3(w)} \quad .
$$

Note that $F(a_1,a_2,a_3)(1,1,1)=|\W(a_1,a_2,a_3)|=(a_1+a_2+a_3)!/(a_1!a_2!a_3!)$.

Given a Laurent polynomial in the variables $q_1,q_2,q_3$, let's call it $P(q_1,q_2,q_3)$, 
let $POS(P)$ be the sum of the monomials all whose
powers are strictly positive, i.e.
$$
POS( \sum_{(i_1,i_2,i_3)} a_{i_1,i_2,i_3} q_1^{i_1}q_2^{i_2}q_3^{i_3} ) \, := \,
\sum_{(i_1,i_2,i_3) ,i_1>0,i_2>0,i_3>0} a_{i_1,i_2,i_3} q_1^{i_1}q_2^{i_2}q_3^{i_3}  \quad .
$$
For example
$$
POS( 5 q_1^{-1} q_2^3 q_3^5 + 4 q_1^{2} q_2^{-3} q_3^5 + 7q_1q_2q_3^2 + 11q_1^2q_3^3 + 2q_1q_2q_3)=  7q_1q_2q_3^2  + 2q_1q_2q_3 \quad .
$$

Defining
$$
G(a_1,a_2,a_3)(q_1,q_2,q_3) \, = \, POS(F(a_1,a_2,a_3)(q_1,q_2,q_3)) \quad ,
$$
the desired number of sucker's bets with three decks of sizes $a_1,a_2,a_3$ and with cards carrying denominations $\{1,2, \dots, a_1+a_2+a_3\}$, is
$$
G(a_1,a_2,a_3)(1,1,1) \quad .
$$

{\bf How to compute $F(a_1,a_2,a_3)$ (and hence  $G(a_1,a_2,a_3)$, and hence  $G(a_1,a_2,a_3)(1,1,1)$ )?}

Every non-empty word $w$ either ends with a `$1$', or with a `$2$', or with a `$3$',
hence the set of words in $1^{a_1}2^{a_2}3^{a_3}$, $\W(a_1,a_2,a_3)$, may be written as
$$
\W(a_1, a_2, a_3) \, = \, \W(a_1-1, a_2, a_3) \, 1 \, \cup \,
\W(a_1, a_2-1, a_3) \, 2 \, \cup \,
\W(a_1, a_2, a_3-1) \, 3 \quad .
$$
Of course $\W(a_1,a_2,a_3)$ is the empty set if any of $a_1,a_2,a_3$ is negative, and $\W(0,0,0)$ consists of
one element, the {\bf empty word}.

Any word $w \in \W(a_1-1,a_2,a_3)$ (given a word $w$ and a letter $l$, $wl$ is the word obtained from $w$ by appending  $l$ to it)
$$
s_1(w1)=s_1(w)+a_2 \quad, \quad
s_2(w1)=s_2(w) \quad, \quad
s_3(w1)=s_3(w)-a_3 \quad .
$$

Similarly, for any word $w \in \W(a_1,a_2-1,a_3)$,
$$
s_1(w2)=s_1(w)-a_1 \quad, \quad
s_2(w2)=s_2(w)+ a_3 \quad, \quad
s_3(w2)=s_3(w) \quad ;
$$
and, for any word $w \in \W(a_1,a_2,a_3-1)$,
$$
s_1(w3)=s_1(w) \quad, \quad
s_2(w3)=s_2(w)- a_2 \quad, \quad
s_3(w3)=s_3(w)+a_1 \quad .
$$
It follows that, 

$\bullet$  For $w \in \W(a_1-1,a_2,a_3)$ we have
$$
weight(w1)= weight(w) \cdot q_1^{a_2} q_3^{-a_3} \quad ;
$$

$\bullet$  For $w \in \W(a_1,a_2-1,a_3)$ we have
$$
weight(w2)= weight(w) \cdot q_1^{-a_1} q_2^{a_3} \quad ;
$$

$\bullet$  For $w \in \W(a_1,a_2,a_3-1)$ we have
$$
weight(w3)= weight(w) \cdot q_2^{-a_2} q_3^{a_1} \quad .
$$

It follows that the Laurent polynomials $F(a_1,a_2,a_3)=F(a_1,a_2,a_3)(q_1,q_2,q_3)$ satisfy the {\bf recurrence relation}
$$
F(a_1,a_2,a_3) \, = \,
q_1^{a_2} q_3^{-a_3} F(a_1-1,a_2,a_3) +
q_1^{-a_1} q_2^{a_3} F(a_1,a_2-1,a_3) +
q_2^{-a_2} q_3^{a_1} F(a_1,a_2,a_3-1) \quad ,
\eqno(Qrecurrence)
$$
subject to the {\it boundary conditions} $F(0,0,0)=1$ and $F(a_1,a_2,a_3)=0$ if $a_1<0$ or $a_2<0$ or $a_3<0$.

This recurrence was programmed in Maple, and from this we deduced the positive parts, $G(a_1,a_2,a_3)$, and plugging-in
$q_1=1,q_2=1,q_3=1$ we obtained

{\bf Important Fact}: The first $12$ terms of the sequence 

`number of sucker's bets' with a set of three decks, each with $n$ cards, where the $3n$ cards have all different numbers 
(labeled $1, \dots, 3n$), starting at $n=1$, are
$$
0, 0, 15, 39, 5196, 32115, 2093199, 19618353, 960165789, 11272949151, 479538890271, 6504453085104 \quad .
$$

Dividing by $3$ to account for trivial cyclic symmetry, the reduced numbers are
$$
0, 0, 5, 13, 1732, 10705, 697733, 6539451, 320055263, 3757649717, 159846296757,    2168151028368 \quad .
$$

It follows, that the sequence of probabilities for a random set of $3$ decks of cards each with $n$ cards to
be a sucker's bet set, for $n$ from 1 to 12 are
$$
0., 0., 0.008928571429, 0.001125541126, 0.006866149723, 0.001872252397,
    0.005245153668, 
$$
$$
0.002072614083, 0.004213592531, 0.002030797274,
    0.003512410777, 0.001921704153 \quad .
$$

{\bf From Counting to Listing using Symbol Crunching}

Suppose that we actually want to see {\it all} the possible sets of three decks that are sucker's bets?
It is probably beyond the scope of computer-kind to list all $6504453085104$ sucker's bets with $3$ decks
of $12$ cards each, but a minor tweak to the recurrence $(Qrecurrence)$ will enable us to do it for smaller sizes
(and in principle, for all sizes).

First we need to tweak the definition of {\it weight}

$$
weightX(w):=q_1^{s_1(w)} q_2^{s_2(w)} q_3^{s_3(w)} \cdot \prod_{i=1}^{|w|} x[i, w_i]\quad .
\eqno(WeightX)
$$

This keeps track of the individuality of the word. Next we define a Laurent polynomial in $q_1,q_2,q_3$ {\bf and}
polynomial in the indeterminates $x[i,j]$ ($1 \leq i \leq |w|$, $j\in \{1,2,3\}$)
$$
F_X(a_1,a_2,a_3)(q_1,q_2,q_3,\{x[i,j]\}) \, := \, \sum_{w \in \W(a_1,a_2,a_3))} weightX(w) \quad .
$$

The same argument that lead to  $(Qrecurrence)$ leads to
$$
F_X \, (a_1,a_2,a_3) \, = \,
$$
$$
x[a_1+a_2+a_3, 1]\, q_1^{a_2} q_3^{-a_3} F_X \, (a_1-1,a_2,a_3) 
$$
$$
+ \, x[a_1+a_2+a_3, 2]\,  q_1^{-a_1} q_2^{a_3} F_X \, (a_1,a_2-1,a_3) 
$$
$$
+\, x[a_1+a_2+a_3, 3]\,  q_2^{-a_2} q_3^{a_1} F_X \, (a_1,a_2,a_3-1) \quad ,
\eqno(Xrecurrence)
$$
subject to the {\it boundary conditions} $F(0,0,0)=1$ and $F(a_1,a_2,a_3)=0$ if $a_1<0$ or $a_2<0$ or $a_3<0$.

We are only interested in the part where all the powers of the variables $q_1,q_2,q_3$ are positive, so let's apply $POS$:
$$
G_X(a_1,a_2,a_3)(q_1,q_2,q_3, \{x[i,j]\}) \, = \, POS(F_X(a_1,a_2,a_3)(q_1,q_2,q_3, \{x[i,j]\})) \quad .
$$
Now each monomial corresponds to a `sucker's bets word'. The computer automatically transcribes each
such monomial to a set of three decks of cards.

See the output file

$\bullet$ {\tt http://sites.math.rutgers.edu/\~{}zeilberg/tokhniot/oSuckerBets3.txt} \quad ,

for the five such sets (up to trivial cyclic symmetry) with 3 decks, each with $3$ cards.

$\bullet$ {\tt http://sites.math.rutgers.edu/\~{}zeilberg/tokhniot/oSuckerBets4.txt} \quad ,

for the thirteen such sets (up to trivial cyclic symmetry) with 3 decks, each with $4$ cards.

$\bullet$ {\tt http://sites.math.rutgers.edu/\~{}zeilberg/tokhniot/oSuckerBets5.txt} \quad ,

for the  $1732$ such sets (up to trivial cyclic symmetry) with 3 decks, each  with $5$ cards.

{\bf From Words to Lattice Paths}

From words, in turn, we can get lattice walks, with unit positive steps, from the origin to the point $(a_1, \dots, a_k)$,
where $w_i=j$ corresponds to the $i$-th step being parallel to  the $x_j$-axis.
For example, the above (`magic-square' word), $132321213$ . corresponds to the 3D lattice walk
$$
(0,0,0) \rightarrow (1,0,0) \rightarrow (1,0,1) 
\rightarrow (1,1,1) \rightarrow (1,1,2) 
\rightarrow (1,2,2) \rightarrow (2,2,2) 
\rightarrow (2,3,2) \rightarrow (3,3,2) 
\rightarrow (3,3,3)  \quad .
$$

From now on the number of decks, $k$, is no longer $3$.

The advantage of the lattice path representation is that it allows us to construct sucker's bets where
repeated denominations are allowed. If  one wants to avoid the possibility of a tie
(like in {\it paper-scissors-stone}) then the allowed steps are always parallel to one of the $k$ axes,
but not necessarily unit steps.
If one does not mind ties, then we can also have `diagonal' steps, where more than one coordinate changes.

For example, the Efron set of four dice mentioned at the very beginning of this treatise corresponds to the 
following lattice path in the 4D hyper-cubic lattice, 
with seven steps
$$
(0,0,0,0) \rightarrow (2,0,0,0) \rightarrow (2,0,0,3) \rightarrow (2,0,4,3) 
\rightarrow (2,6,4,3) \rightarrow (6,6,4,3) \rightarrow (6,6,4,6) \rightarrow (6,6,6,6) \quad .
$$

One can easily define analogs of $s_i(w)$ for these more general walks (bets), and establish the
generalizations of $(Qrecurrence)$ and $(Xrecurrence)$. See the Maple package {\tt SuckerBets.txt}
available directly from

{\tt http://sites.math.rutgers.edu/\~{}zeilberg/tokhniot/SuckerBets.txt}  \quad ,

or via the front of this treatise

{\tt http://sites.math.rutgers.edu/\~{}zeilberg/mamarim/mamarimhtml/suckerbets.html} \quad ,

where there are links to plenty of input and output files.

The Efron set of dice, mentioned above, is a bit wasteful, and there exists a  tie-less four-dice sucker's bet set
only using {\bf six} different denominations, see

{\tt http://sites.math.rutgers.edu/\~{}zeilberg/tokhniot/oSuckerBets0a.txt} \quad.

Here it is:
$$
A=[1,1,5,5,5,5] \quad , \quad
B=[4,4,4,4,4,4] \quad , \quad
C=[3,3,3,3,3,3] \quad , \quad
D=[2,2,2,2,6,6] \quad .
$$

This is the {\bf unique} such set, up to trivial cyclic symmetry.

Efron's set is one of $38$ different (tie-less, and up to trivial cyclic symmetry) sets with {\bf seven} different denominations, viewable from

{\tt http://sites.math.rutgers.edu/\~{}zeilberg/tokhniot/oSuckerBets0b.txt} \quad .

There are $755$ different (tie-less, and up to trivial cyclic symmetry) sets with {\it eight} different denominations, viewable from

{\tt http://sites.math.rutgers.edu/\~{}zeilberg/tokhniot/oSuckerBets0c.txt} \quad .

For many more examples, where the number of cards (or faces) in each deck (or die) are not necessarily the same, see
the other outputs files.

{\bf Statistical Analysis}

Now we are back to $k=3$. What we do here for three-deck sets is generalizable to general $k \geq 3$,
but things are already interesting (and complex enough) for the case $k=3$.

The three individual {\it word-statistics} $s_1(w), s_2(w), s_3(w)$ are closely related, and in fact,
trivially equivalent to, the so-called `Number of inversions' introduced  in  the 19th century by
Eugen Netto (of the {\it Lehrbuch} fame), and rediscovered, at the mid 20th century, by non-parametric statisticians
H. B. Mann and D. R. Whitney ([MW]). It is well known (see, e.g. [CJZ]), that each of $s_1(w),s_2(w), s_3(w)$
is, {\it individually}, {\bf asymptotically normal} , but of course, they are far from being independent.

We discovered that the {\it limiting scaled} tri-variate distribution of the triple (discrete) random variable $(s_1(w),s_2(w),s_3(w))$
defined on the set of words on $1^n2^n3^n$, namely $\W(n,n,n)$,  as $n \rightarrow \infty$,  converges, {\bf in distribution}, to the
{\bf limit} as $c \rightarrow 1^{-}$ of the trivariate continuous random variable whose {\bf  joint density} function is
$$
f(x,y,z; c)\, := \, \frac{exp(-x^2/-y^2/2-z^2/2 \,- \, c\, (xy+xz+yz) )}{N(c)} \quad ,
$$
where $N(c)$ is the normalization factor that would make $\int_{R^3} f(x,y,z;c)=1$, namely
$$
N(c) := \frac{(2\pi)^{3/2}}{(1-c)\sqrt{1+2c}} \quad .
$$
Note that this `blows up' at $c=1$.

Officially this is still a conjecture, but we are {\bf absolutely sure}  that it is correct, since we proved,
{\bf rigorously}, that all the scaled mixed moments, $M(i_1,i_2,i_2)$,  of the triple of random variables $(s_1,s_2,s_3)$, converge as
$n \rightarrow \infty$ to the limit of the corresponding scaled mixed moments of $f(x,y,z;c)$, as $c \rightarrow 1^{-}$,
for all $1 \leq i_1,i_2,i_3 \leq 5$, and with more computing power, one can easily go further.

We are offering to donate $100$ dollars  to the OEIS, in honor of the first prover of our conjecture.
(Not because we have any doubts about its truth, or that we particularly care about so-called rigorous proofs,
but because we love the OEIS!)

{\bf Explicit (Rigorously Proved!) Expressions for Mixed Moments of $(s_1,s_2,s_3)$}

See: {\tt http://sites.math.rutgers.edu/\~{}zeilberg/tokhniot/oSuckerBetsAnalysis4LP.txt} \quad .

Let's just summarize a few highlights.

The variance  of each of $s_1,s_2,s_3$  on $\W(n,n,n)$ is ${n}^{2} \left( 2\,n+1 \right)/3$, and
the kurtosis is 
$$
\frac{3(10n^2-n-4)}{5n(2n+1)} \quad .
$$
Note that it tends to $3$, as $n$ goes to infinity, as it should, since it is {\it asymptotically normal}.
This is all old stuff (see [CJZ]), as well as any of the moments of each of the single random variables $s_1,s_2,s_3$.

The {\it covariance} between any pair of $s_1,s_2,s_3$ is very simple, it is
$$
-\frac{n^3}{3} \quad,
$$
hence the  {\it correlation is}
$$
-\frac{n}{2n+1} \quad ,
$$
that converges to $-\frac{1}{2}$. 

Explicit expressions for all the mixed moments can be found in the above output file.

Let's just write down, in humanese, the $(4,5,5)$ mixed moment, that happens to be a polynomial in $n$ of degree $21$.
$$
{\frac{1}{2837835}}\,{n}^{3} ( 39239200\,{n}^{18}+66146080\,{n}^{17}-816055240\,{n}^{16}
$$
$$
+1114633520\,{n}^{15}+3208398492\,{n}^{14}-13589761044\,{n}^{13}+25028291837\,{n}^{12}-38043392560\,{n}^{11}+62580129596\,{n}^{10}
$$
$$
-103184180072\,{n}^{9}+157753326632\,{n}^{8}-224678523360\,{n}^{7}+293133737664\,{n}^{6}-336053442624\,{n}^{5}
$$
$$
+322828696448\,{n}^{4}-243844376832\,{n}^{3}+132045454336\,{n}^{2}-44452356096\,n+6864979968) \quad .
$$

Here are all the scaled limits, $S(i_1,i_2,i_3)$, for $0 \leq i_1,i_2,i_3 \leq 5$. By symmetry
we only need to list the values for $i_1 \leq i_2 \leq i_3$. Note that $S(i_1,i_2,i_3)=0$ when $i_1+i_2+i_3$ is odd.
$$
S(0,0,0) = 1 \quad , \quad
S(0,0,2) = 1 \quad , \quad
S(0,0,4) = 3 \quad , \quad
S(0,1,1) = -1/2 \quad , \quad
$$
$$
S(0,1,3) = -3/2 \quad , \quad
S(0,1,5) = -15/2 \quad , \quad
S(0,2,2) = 3/2 \quad , \quad
S(0,2,4) = 6 \quad , \quad
$$
$$
S(0,3,3) = -21/4 \quad , \quad
S(0,3,5) = -30 \quad , \quad
S(0,4,4) = 57/2 \quad , \quad
S(0,5,5) = -765/4 \quad , \quad
$$
$$
S(1,1,2) = 0 \quad , \quad
S(1,1,4) = 3/2 \quad , \quad
S(1,2,3) = -3/4 \quad , \quad
S(1,2,5) = -15/2 \quad , \quad
$$
$$
S(1,3,4) = 3/2 \quad , \quad
S(1,4,5) = -45/4 \quad , \quad
S(2,2,2) = 3/2 \quad , \quad
S(2,2,4) = 6 \quad , \quad
$$
$$
S(2,3,3) = -3 \quad , \quad
S(2,3,5) = -45/2 \quad , \quad
S(2,4,4) = 45/2 \quad , \quad
S(2,5,5) = -135 \quad , \quad
$$
$$
S(3,3,4) = 0 \quad , \quad
S(3,4,5) = -135/4 \quad , \quad
S(4,4,4) = 135/2 \quad , \quad
S(4,5,5) = -945/4 \quad  .
$$

They do indeed  coincide with the llimits as $c \rightarrow 1^{-}$ of the scaled moments of $f(x,y,z;c)$.
The mixed moments of the latter can be computed as far as desired thanks to linear recurrences
of order $4$ found via the Apagodu-Zeilberger multi-variable Almkvist-Zeilberger algorithm ([AZ]).
When one plugs-in $c=1$ you still get complicated recurrences, but the {\bf diagonal}, the mixed
moments $(2n,2n,2n)$ are given by a very nice closed form
$$
S(2n,2n,2n) \, = \, {\frac { \left( 3\,n \right) !\, \left( 2\,n \right) !}{{8}^{n} \left( n! \right) ^{2}}} \quad .
$$

{\bf How To Find Rigorously Proved Polynomial Expressions for the Mixed Moments of $(s_1,s_2,s_3)$}?

We use the approach  of [BZ] and [Z].  Suppose that we are interested in an explicit expression for a certain
specific mixed moment $M(i,j,k)$ of the triple $(s_1,s_2,s_3)$ defined on $\W(n,n,n)$. 
We know, {\it a priori},  that it is a certain {\it polynomial} in $n$, and we can easily bounds its degree, let's call it $d$.
Hence it suffices to find $d+1$ terms. But $M(i,j,k)$ at $(n,n,n)$ is
$$
(q_1 \frac{\partial}{\partial q_1})^i
(q_2 \frac{\partial}{\partial q_2})^j  
(q_3 \frac{\partial}{\partial q_3})^k
F(n,n,n)(q_1,q_2,q_3) \Bigl \vert_{q_1=1,q_2=1,q_3=1} \quad .
$$
We use $(Qrecurrence)$ to compute enough terms $F(n,n,n)$ for $n\leq d+1$, and then `fit the data'.

The drawback of the above approach is that it is inefficient. If we are only interested in the first few mixed moments,
then we can do a tri-variate Taylor expansion around $(q_1,q_2,q_3)=(1,1,1)$ and truncate anything beyond
our horizon. This is accomplished by using $(Qrecurrence)$, with $q_1=1+p_1, q_2=1+p_2, q_3=1+p_3$ and
doing the truncated version. This gives a quicker way to find the {\it mixed factorial moments} from which the mixed moments can
be easily gotten.
Full details can be gotten by reading the Maple source code available
from the front of this article mentioned above.

{\bf Conclusion} 

As with most of our research, the {\it methodology} is at least as interesting as
the actual results. The present project is a {\bf case study} of using {\it symbol crunching} and {\it experimental mathematics}
both to generate interesting combinatorial objects (in this case {\it sucker's bets}), via {\it symbolic dynamical programming},
and to do {\it symbolic statistical analysis}  about intriguing combinatorial random variables, that
are very negatively correlated, and for which there is a beautiful {\it asymptotic} limiting trivariate distribution, namely
the limit, as $c$ goes to $1$ from the below, of the one whose joint density function is
$$
f(x,y,z; c)\, := \, \frac{exp(-x^2/-y^2/2-z^2/2 \,- \, c\, (xy+xz+yz) )}{ \frac{(2\pi)^{3/2}}{(1-c)\sqrt{1+2c}}} \quad .
$$

{\bf References}

[AZ] Moa Apagodu and Doron Zeilberger, {\it Multi-Variable Zeilberger and Almkvist-Zeilberger Algorithms and the Sharpening of Wilf-Zeilberger Theory},
Adv. Appl. Math. {\bf 37}(2006), (Special issue in honor of Amitai Regev), 139-152. \hfill\break
{\tt http://sites.math.rutgers.edu/\~{}zeilberg/mamarim/mamarimhtml/multiZ.html}

[BZ] Andrew Baxter and Doron Zeilberger,
{\it The Number of Inversions and the Major Index of Permutations are Asymptotically Joint-Independently Normal (2nd ed.)},
Personal Journal of Shalosh B. Ekhad and Doron Zeilberger, Feb. 4, 2011. \hfill\break
{\tt http://sites.math.rutgers.edu/\~{}zeilberg/mamarim/mamarimhtml/invmaj.html} \quad .

[CJZ]  E. Rodney Canfield, Svante Janson, and Doron Zeilberger,
{\it The Mahonian Probability Distribution on Words is Asymptotically Normal},
Advances in Applied Mathematics {\bf 46}(2011), 109-124. Available from: \hfill\break
{\tt http://sites.math.rutgers.edu/\~{}zeilberg/mamarim/mamarimhtml/mahon.html}  \quad , \hfill\break
erratum: \hfill\break
{\tt http://sites.math.rutgers.edu/\~{}zeilberg/mamarim/mamarimPDF/mahonerratum.pdf} \quad .

[G] Martin Gardner, 
{\it ``The Colossal Book of Mathematics: Classic Puzzles, Paradoxes, and Problems: 
Number Theory, Algebra, Geometry, Probability, Topology, Game Theory, Infinity, and Other Topics of Recreational Mathematics''} (1st ed.),
New York: W. W. Norton \& Company..

[MW] H. B. Mann and D. R. Whitney, {\it On a test whether one of two random
variables is stochastically larger than the other}. Annals of Mathematical
Statistics {\bf 18}(1947), 50-60.

[NZ] John Noonan and Doron Zeilberger, {\it ``The Goulden-Jackson Cluster Method: Extensions, Applications, and Implementations"},
J. Difference Eq. Appl. {\bf 5} (1999), 355-377. \hfill\break
{\tt http://sites.math.rutgers.edu/\~{}zeilberg/mamarim/mamarimhtml/gj.html} \quad .

[P] W. Penney, {\it Problem 95: Penney ante},  J. Recreational Mathematics {\bf 7}(1974), 321.

[W] The Wikipedia Foundation, {\it ``Nontransitive Dice''}, \hfill\break
{\tt https://en.wikipedia.org/wiki/Nontransitive\_dice} \quad .

[Z] Doron Zeilberger, {\it Doron Gepner's Statistics on Words in $\{1,2,3\}^{*}$ is (Most Probably) Asymptotically Logistic}.
Personal Journal of Shalosh B. Ekhad and Doron Zeilberger, March 31, 2016. \hfill\break
{\tt http://sites.math.rutgers.edu/\~{}zeilberg/mamarim/mamarimhtml/gepner.html} \quad .

\bigskip
\bigskip
\hrule
\bigskip
Doron Zeilberger, Department of Mathematics, Rutgers University (New Brunswick), Hill Center-Busch Campus, 110 Frelinghuysen
Rd., Piscataway, NJ 08854-8019, USA. \hfill \break
DoronZeil at gmail dot com  \quad ;  \quad {\tt http://sites.math.rutgers.edu/\~{}zeilberg/} \quad .
\bigskip
\hrule
\bigskip
Shalosh B. Ekhad, c/o D. Zeilberger, Department of Mathematics, Rutgers University (New Brunswick), Hill Center-Busch Campus, 110 Frelinghuysen
Rd., Piscataway, NJ 08854-8019, USA.
\bigskip
\hrule

\bigskip
Exclusively published in The Personal Journal of Shalosh B. Ekhad and Doron Zeilberger  \hfill \break
{ \tt http://sites.math.rutgers.edu/\~{}zeilberg/pj.html} \quad  and {\tt arxiv.org} \quad . 
\bigskip
\hrule
\bigskip
{\bf  Written:  Oct. 27, 2017.}

\end